\documentclass[twoside]{article}
\usepackage{authblk,amsmath,amssymb,mathrsfs,latexsym,amsthm,geometry}
\theoremstyle{definition}
\newtheorem{Def}{Definition}[subsection]
\newtheorem{Ex}[Def]{Example}
\newtheorem{Thm}[Def]{Theorem}
\newtheorem{Lem}[Def]{Lemma}
\newtheorem{Prop}[Def]{Proposition}
\newtheorem{Cor}[Def]{Corollary}
\newtheorem{Rem}[Def]{Remark}
\newtheorem{Main}{Main Theorem}
\newtheorem{M}{Main Theorem}
\newtheorem{Prf}{Proof}
\newtheorem{Th}{Theorem \ref{Id}}

\title{Vanishing theorems of $L^2$-cohomology groups on Hessian manifolds }
\author{Shinya Akagawa}
\affil{Department of Mathematics, Graduate School of Science, Osaka University, Osaka 560, Japan}
\date{}
\begin{document}

\maketitle

\begin{abstract}
We show vanishing theorems of $L^2$-cohomology groups of Kodaira-Nakano type on complete Hessian manifolds. We obtain further vanishing theorems of $L^2$-cohomology groups $L^2H^{p,q}(\Omega)$ on a regular convex cone $\Omega$ with the Cheng-Yau metric for $p>q$.
\end{abstract}

\begin{enumerate}
\item[]\emph{keywords}: Hessian manifolds, Hesse-Einstein, Monge-Amp\`{e}re equation, Laplacians, $L^2$-cohomology groups, regular convex cones.

\item[]\emph{MS classification}: 53C25, 53C55. 
\end{enumerate}

\tableofcontents

\setcounter{section}{-1}

\section{Introduction}

A \emph{flat manifold} $(M,D)$ is a manifold $M$ with a flat affine connection $D$, where an affine connection is said to be \emph{flat} if the torsion and the curvature vanish identically. A flat affine connection $D$ gives an \emph{affine local coordinate system} $\{x^1, \dots , x^n\}$ satisfying 
\[D_{\frac{\partial}{\partial x^i}} \frac{\partial}{\partial x^j} = 0.\]
A Riemannian metric $g$ on a flat manifold $(M,D)$ is said to be a \emph{Hessian metric} if $g$ can be locally expressed in the Hessian form
with respect to an affine coordinate system $\{x^1, \dots , x^n\}$ and a \emph{potential function} $\varphi$, that is,
\[g_{ij}=\frac{\partial^2 \varphi}{\partial x^i \partial x^j}.\]
The triplet $(M,D,g)$ is called a \emph{Hessian manifold}. The Hessian structure $(D,g)$ induces a holomorphic coordinate system $\{ z^1, \dots ,z^n \}$ and a K\"{a}hler metric $g^T$ on $TM$ such that
\[z^i=x^i+ \sqrt{-1} y^i,\]
\[g^T_{i\bar{j}}(z)=g_{ij}(x),\]
where $\{ x^1, \dots ,x^n,y^1, \dots ,y^n\}$ is a local coordinate system on $TM$ induced by the affine coordinate system $\{x^1, \dots, x^n\}$ and fibre coordinates $\{y^1, \dots , y^n\}$. In this sense, Hessian geometry is a real analogue of K\"{a}hler geometry.

A $(p,q)$-form on a flat manifold $(M,D)$ is a smooth section of $\wedge^p T^*M \otimes \wedge^q T^*M$. On the space of $(p,q)$-forms, a flat connection $D$ induces a differential operator $\bar{\partial}$ which is an analogue of the Dolbeault operator. Then the cohomology group $H^{p,q}_{\bar{\partial}}(M)$ is defined with respect to $\bar{\partial}$. On compact Hessian manifolds, Shima proved an analogue of Kodaira-Nakano vanishing theorem for $H^{p,q}_{\bar{\partial}}(M)$ by using the theory of harmonic integrals (c.f. Theorem \ref{C} \cite{Sh}). However, most of important examples of Hessian manifolds such as regular convex domains (c.f. Theorem \ref{M} \cite{CY}) are noncompact. Therefore, we prove an extension of Theorem \ref{C} in the case of complete Hessian manifolds in Section 3.2.

\setcounter{M}{1}

\begin{M}
Let $(M,D,g)$ be an oriented $n$-dimensional complete Hessian manifold and $(F,D^F)$ a flat line bundle over $M$. We denote by $h$ a fiber metric on $F$. Assume that there exists $\varepsilon >0$ such that $B+ \beta= \varepsilon g$ where $B$ and $\beta$ are the second Koszul forms with respect to fiber metric $h$ and Hessian metric $g$ respectively. Then if for $p+q>n$ and all $v \in L ^{p,q}(M,g,F,h)$ such that $ \bar { \partial }v=0$, there exists $u \in L^{p,q-1}(M,g,F,h)$ such that 
\[ \bar { \partial }u=v, \quad  \| u \| \le  \{ \varepsilon(p+q-n)\}^{- \frac {1}{2}} \| v \| .\] 
In particular, we have
 \[L^2H^{p,q}_{\bar{\partial}}(M,g,F,h)=0, \quad \textrm{for} \ p+q>n.\]
\end{M}

Remark that we cannot use the harmonic theory for the proof and we need the method of functional analysis as in the case of complete K\"{a}hler manifolds. To prove Main Theorem \ref{LC}, we introduce the operator $\partial^\prime_F$ (c.f. Definition \ref{Pr}) which is not defined in \cite{Sh} and we obtain the following as an analogue of Kodaira-Nakano identity.

\begin{Th}
Let $(D,g)$ is a Hessian structure. Then we have
\[ \bar { \square }_F= \square ^\prime _F +[ e( \beta +B ), \Lambda ]. \]
\end{Th}

An open convex cone $\Omega$ in $\mathbb{R}^n$ is said to be \emph{regular} if $\Omega$ contains no complete straight lines. We can apply Main Theorem \ref{LC} to regular convex cones with the \emph{Cheng-Yau metric} (c.f. Theorem \ref{M} \cite{CY}). Further, we have stronger vanishing theorems as follows in Section 3.3.

\begin{M}
Let $(\Omega,D,g=Dd\varphi)$ be a regular convex cone in $\mathbb{R}^n$ with the Cheng-Yau metric. Then for $p>q$ and all $v \in L ^{p,q}(M,g)$ such that $ \bar { \partial }v=0$, there exists $u \in L^{p,q-1}(M,g)$ such that 
\[ \bar { \partial }u=v, \quad  \| u \| \le  (p-q)^{- \frac {1}{2}} \| v \| .\]
In particular, we have 
\[L^2H^{p,q}_{\bar{\partial}}(M,g)=0, \quad \textrm{for} \ p>q.\]
\end{M}

In the case of a Hessian manifold $(\mathbb{R}^n,D,g)$ as in Example \ref{E} $(2)$, we have sharp vanishing theorem in Section 3.4.

\begin{M}
For  $p \ge 1$, $q \ge 0$ and $v \in L ^{p,q}(\mathbb{R}_+^n,g)$ such that $\bar { \partial }v=0$, there exists $u \in L^{p,q-1}(\mathbb{R}_+^n,g)$ such that 
\[ \bar { \partial }u=v, \quad  \| u \| \le p^{- \frac {1}{2}} \| v \| .\]
In particular, we have 
\[L^2H^{p,q}_{\bar{\partial}}(\mathbb{R}_+^n,g)=0, \quad \textrm{for} \ p \geq 1 \ \textrm{and} \ q \geq 0.\]
\end{M}

\section{Hessian manifolds} \label{Fu}
\subsection{Hessian manifolds}

\begin{Def}
An affine connection $D$ on a manifold $M$ is said to be \emph{flat} if the torsion tensor $T^D$ and the curvature tensor $R^D$ vanish identically. A manifold $M$ endowed with a flat connection $D$ is called a \emph{flat manifold}, which is denoted by $(M,D)$.
\end{Def}

\begin{Prop} \label{af}
\cite{Sh} Suppose that $(M,D)$ is a flat manifold. Then there exists a local coordinate system $\{x^1, \dots ,x^n\}$ on $M$ such that ${D_{\frac{\partial}{\partial x^i}}}\frac{\partial}{\partial x^j}=0$. The changes between such a local coordinate system are affine transformations.
\end{Prop}

\begin{Def}
The local coordinate system in Proposition \ref{af} is called an \emph{affine coordinate system} with respect to $D$.
\end{Def}

In this paper, every local coordinate system on flat manifolds is given as an affine coordinate system.

\begin{Def}
A Riemannian metric $g$ on a flat manifold $(M,D)$ is said to be a \emph{Hessian metric} if $g$ is locally expressed by
\[g=Dd \varphi ,\]
that is,
\[g_{ij}=\frac{\partial^2 \varphi}{\partial x^i \partial x^j}.\]
Then the pair $(D,g)$ is called a \emph{Hessian structure} on $M$, and $\varphi$ is said to be a potential of $(D,g)$. A manifold $M$ with a Hessian structure $(D,g)$ is called a \emph{Hessian manifold}, which is denoted by $(M,D,g)$.
\end{Def}

Let $(M,D)$ be a flat manifold and $TM$ the tangent bundle over $M$. We denote by $\{x^1, \dots ,x^n,y^1, \dots ,y^n\}$ a local coordinate system on $TM$ induced by an affine coordinate system $\{x^1, \dots ,x^n\}$ on $M$ and fibre coordinates $\{y^1, \dots , y^n\}$. Then a holomorphic coordinate system$ \{z^1, \dots ,z^n\}$ on $TM$ is given by
\[z^i=x^i + \sqrt{-1} y^i.\]
For a Riemannian metric $g$ on $M$ we define a Hermitian metric $g^T$ on $TM$ by
\[g^T = \sum_{i,j}g_{ij}dz^i \otimes d\bar{z}^j.\]

\begin{Prop} \cite{Sh}
Let $(M,D)$ be a flat manifold and $g$ a Riemannian metric on $M$. Then the following conditions are equivalent.
\item[$(1)$] $g$ is a Hessian metric.
\item[$(2)$] $g^T$ is a K\"{a}hler metric.
\end{Prop}

\begin{Ex} \label{E}
\item[$(1)$] Let $(D,g)$ be a pair consisting of the standard affine connection $D$ and a Euclidean metric on $\mathbb{R}^n$. Then $(D,g)$ is a Hessian structure. Indeed, if we set $\displaystyle \varphi (x) = \frac{1}{2} \sum_j (x^j)^2$, we have
\[\displaystyle \frac{\partial ^2 \varphi}{\partial x^i \partial x^j} = \delta_{ij} = g_{ij},\]
where $\delta_{ij}$ is the Kronecker delta, that is,
\[\delta_{ij} = \begin{cases} 1 & (i=j) \\ 0 & (i \neq j). \end{cases}\]
Moreover, the K\"{a}hler metric $g^T$ on $T \mathbb{R}^n \simeq \mathbb{C}^n$ is also a Euclidean metric.
\item[$(2)$] We set $\mathbb{R}_+ = (0,\infty)$. Let $D$ be the standard affine connection, that is, the restriction of a the standard affine connection on $\mathbb{R}^n$ to $\mathbb{R}_+^n$. We define a Riemannian metric $g$ on $\mathbb{R}_+ ^n$ by 
\[ g_{ij}(x)= \frac{\delta_{ij}}{(x^j)^2}.\]
Then $(D,g)$ is a Hessian structure. Indeed, if we set $\varphi (x)=-\log (x^1 \cdots x^n)$, we have
\[ \frac{\partial ^2 \varphi}{\partial x^i \partial x^j} = g_{ij}.\]
In addition, the K\"{a}hler metric $g^T$ on $T \mathbb{R}_+ \simeq \mathbb{R}_+ \oplus \sqrt{-1} \, \mathbb{R}$ is the Poincar\'e metric.
\end{Ex}

\begin{Def} \label{na}
Let $M$ be a manifold and $D$ a torsion-free affine connection on $M$. We denote by $g$ a Riemannian metric on $M$, and by $\nabla$ the Levi-Civita connection of $g$. We define a \emph{difference tensor} $\gamma$ of $\nabla$ and $D$ by
\[\gamma = \nabla - D.\]
\end{Def}

We denote by $\mathscr{X}(M)$ the space of vector fields on $M$. Since $\nabla$ and $D$ are torsion-free, it follows that for $X,Y \in \mathscr{X}(M)$
\[\gamma_XY=\gamma_YX.\]
It should be remarked that the components ${\gamma^i}_{jk}$ of $\gamma$ with respect to affine coordinate systems coincide with the Christoffel symbols of $\nabla$.

\begin{Prop} \label{H}
\cite{Sh} Let $(M,D)$ be a flat manifold and $g$ a Riemannian manifold on $M$. Then the following conditions are equivalent.
\item[$(1)$] $(D,g)$ is a Hessian structure.
\item[$(2)$]
$\displaystyle (D_Xg)(Y,Z)=(D_Yg)(X,Z) , \quad X,Y,Z \in \mathscr{X}(M) \quad
 \Bigl(\Leftrightarrow \frac{\partial g_{jk}}{\partial x^i} = \frac{\partial g_{ik}}{\partial x^j} \Bigr)$.
\item[$(3)$]
$\displaystyle g(\gamma_XY,Z)=g(Y,\gamma_XZ) , \quad X,Y,Z \in \mathscr{X}(M) \quad
 (\Leftrightarrow\gamma_{ijk}=\gamma_{jik} )$.
\item[$(4)$]
$\displaystyle (D_Xg)(Y,Z)=2g(\gamma_XY,Z) , \quad X,Y,Z \in \mathscr{X}(M) \quad
 \Bigl(\Leftrightarrow \frac{\partial g_{ij}}{\partial x^k} = 2\gamma_{ijk}\Bigr)$.
\end{Prop}

\begin{Def} \label{Du}
Let $M$ be a manifold and $D$ a torsion-free affine connection on $M$. We denote by $g$ a Riemannian metric on $M$. We define another affine connection $D^*$ on $M$ as follows:
\[Xg(Y,Z)=g(D_XY,Z)+g(Y,D^*_XZ), \quad X,Y,Z \in \mathscr{X}(M)\]
We call $D^*$ the \emph{dual connection} of $D$ with respect to $g$.
\end{Def}

\begin{Prop} \label{S}
\cite{Sh} Let $(M,g)$ be a Riemannian manifold and $D$ a torsion-free affine connection on $M$. We denote by $D^*$ the dual connection of $D$ with respect to $g$. Let $\gamma$ be the  difference tensor of $\nabla$ and $D$. Then the following conditions are equivalent.
\item[$(1)$] $D^*$ is torsion-free．
\item[$(2)$] $(D_Xg)(Y,Z)=(D_Yg)(X,Z), \quad X,Y,Z \in \mathscr{X}(M)$.
\item[$(3)$] $g(\gamma_XY,Z)=g(Y,\gamma_XZ) , \quad X,Y,Z \in \mathscr{X}(M)$.
\item[$(4)$] $(D_Xg)(Y,Z)=2g(\gamma_XY,Z) , \quad X,Y,Z \in \mathscr{X}(M)$.
\item[$(5)$] $D+D^*=2\nabla$.
\end{Prop}

\subsection{Koszul forms on flat manifolds}

\begin{Def} \label{DKos}
Let $(M,D)$ be a flat manifold and $g$ a Riemannian metric on $M$. We define a $d$-closed $1$-form $ \alpha$ and a symmetric bilinear form $ \beta$ by
\[\alpha =\frac {1}{2} d \log \det [g_{ij}]  , \quad \beta = D \alpha .\]
Remark that since the changes between affine coordinate systems are affine transformations, $\alpha$ and $\beta$ are globally well-defined. We call $\alpha$ and $\beta$ the \emph{first Koszul form} and the \emph{second Koszul form} for $(D,g)$, respectively.  
\end{Def}

\begin{Prop} \label{Kos}
\cite{Sh} Let $(M,D,g)$ be a Hessian manifold. Then we have the following equations.
\[\alpha_i:=\alpha(\frac{\partial}{\partial x^i})=\sum_r {\gamma^r}_{ri}, \quad
\beta_{ij}:=\beta(\frac{\partial}{\partial x^i},\frac{\partial}{\partial x^j})=\sum_r \frac{\partial {\gamma^r}_{ri}}{\partial x^j}.\]
\end{Prop}

\begin{Def}
Let $(M,D,g)$ be a Hessian manifold. If there exists $\lambda \in \mathbb{R}$ such that $\beta=\lambda g$, we call $g$ a \emph{Hesse-Einstein metric}. 
\end{Def}

It should be remarked that a Hessian metric $g$ on $M$ is a Hesse-Einstein metric if and only if a K\"{a}hler metric $g^T$ on $TM$ is a K\"{a}hler-Einstein metric \cite{Sh}.

A convex domain in $\mathbb{R}^n$ which contains no full straight lines is called a \emph{regular convex domain}. By the following theorem, on a regular convex domain there exists a complete Hesse-Einstein metric $g$ which satisfies $g=\beta$. It is called the \emph{Cheng-Yau metric}.

\begin{Thm} \label{M}
\cite{CY} On a regular convex domain $\Omega \in \mathbb{R}^n$, there exists a unique convex function $\varphi$ such that
\[ 
\begin{cases} 
\det \Bigl[\frac{\partial^2 \varphi}{\partial x^i \partial x^j} \Bigr] =e^{2\varphi} \\
\varphi(x) \rightarrow \infty & (x \rightarrow \partial \Omega) .
\end{cases}
\]
In addition, the Hessian metric $g=Dd \varphi$ is complete, where $D$ is the standard affine connection on $\Omega$.
\end{Thm}

\begin{Cor} \label{A}
The Cheng-Yau metric $g$ defined by Theorem \ref{M} is invariant under affine automorphisms of $\Omega$, where an affine automorphism of $\Omega$ is restriction of an affine transformation $A: \mathbb{R}^n \to \mathbb{R}^n$ to $\Omega$ which satisfies $A\Omega=\Omega$.
\end{Cor}

\begin{Prf}
An affine transformation $A$ is denoted by
\[Ax=((Ax)^1, \dots ,(Ax)^n), \quad (Ax)^i= \sum_j a^i_j x_j + b^i.\]
We define a function $\tilde{\varphi}$ on $\Omega$ by
\[\tilde{\varphi}(x) = \varphi(Ax) + \log | \det [a^i_j] | .\]
Then we have
\[\tilde{\varphi}(x) \rightarrow \infty \quad (x \rightarrow \infty).\]
Moreover we obtain
\[\frac{\partial^2 \tilde{\varphi}}{\partial x^i \partial x^j}(x)
= \sum_{k,l} a^k_i a^l_j \frac{\partial^2 \varphi}{\partial x^k \partial x^l}(Ax).\]
Hence $\tilde{\varphi}$ is a convex function. Furthermore, it follows that
\[\det \Bigl[\frac{\partial^2 \tilde{\varphi}}{\partial x^i \partial x^j}(x)\Bigr]
=| \det [a^i_j] |^2 \det \Bigl[\frac{\partial^2 \varphi}{\partial x^i \partial x^j}(Ax)\Bigr]
=e^{2(\varphi(Ax) + \log | \det [a^i_j] |)} =e^{2\tilde{\varphi}(x)}.\]
Therefore $\tilde{\varphi}$ is also a convex function which satisfies the condition of Theorem \ref{M}. From the uniqueness of the solution we have $\tilde{\varphi}=\varphi$, that is,
\[\varphi(x)=\varphi(Ax) + \log | \det [a^i_j] |. \]
Hence we have 
\[ g_{ij}(x)=\sum_{k,l} a^k_i a^l_j g_{kl}(Ax).\]
This implies that $g$ is invariant under affine automorphisms.
\begin{flushright}
$\Box$
\end{flushright}
\end{Prf}

\begin{Ex}
Let $(\mathbb{R}_+^n,D,g=Dd\varphi)$ be the same as in Example \ref{E} $(2)$. Then $\varphi(x)=-\log (x^1 \cdots x^n)$ satisfies the condition of Theorem \ref{M}.
\end{Ex}

\section{$(p,q)$-forms on flat manifolds} \label{PQ}

Hereafter, we assume that $(M,D)$ is an oriented flat manifold and $g$ is a Riemannian metric on $M$. In addition, let $F$ be a real line bundle over $M$ endowed with a flat connection $D^F$ and a fiber metric $h$. Moreover, we denote by $\{s\}$ a local frame field on $F$ such that $D^F s=0$.

\subsection{$(p,q)$-forms and fundamental operators} 

\begin{Def}
We denote by $A^{p,q}(M)$ the space of smooth sections of $\wedge ^p T^* M \otimes \wedge ^q T^* M$. An element in $A^{p,q}(M)$ is called a $(p,q)$-form. For a $p$-form $\omega$ and a $q$-form $\eta$, $\omega \otimes \eta \in A^{p,q}(M)$ is denoted by $\omega \otimes \bar{\eta}$.
\end{Def}

Using an affine coordinate system a $(p,q)$-form $\omega$ is expressed by
\[\omega=\sum_{I_p ,J_q}\omega_{I_p J_q}dx^{I_p}\otimes \overline{dx^{J_q}},\]
where
\[I_p=(i_1, \dots ,i_p), \quad 1 \le i_1< \dots <i_p \le n, \quad J_q=(j_1, \dots ,j_q), \quad 1 \le j_1< \dots <j_q \le n,\]
\[ dx^{I_p}=dx^{i_1}\wedge \dots \wedge dx^{i_p}, \quad dx^{J_q}=dx^{j_1}\wedge \dots \wedge dx^{j_q}.\]

\begin{Ex}
A Riemannian metric $g$ and the second Koszul form $\beta$ (Definition \ref{DKos}) are regarded as $(1,1)$-forms; 
\[g=\sum_{i,j}g_{ij} \, dx^i \otimes \overline{dx^j}, \quad \beta=\sum_{i,j}\beta_{ij} \, dx^i \otimes \overline{dx^j}.\]
\end{Ex}

\begin{Def}
We define the exterior product of $\omega \in A^{p,q}(M)$ and $\eta \in A^{r,s}(M)$ by
\[\omega \wedge \eta =\sum_{I_p ,J_q ,K_r ,L_s}\omega_{I_p J_q}\eta_{K_r L_s}dx^{I_p} \wedge dx^{K_r} \otimes \overline{dx^{J_q} \wedge dx^{L_s}},\]
where $\displaystyle \omega=\sum_{I_p ,J_q}\omega_{I_p J_q}dx^{I_p}\otimes \overline{dx^{J_q}}$ and $\displaystyle \eta=\sum_{K_r ,L_s}\eta_{K_r L_s}dx^{K_r}\otimes \overline{dx^{L_s}}$.
\end{Def}

\begin{Def}
For $\omega \in A^{r,s}(M)$ we define an exterior product operator $e(\omega) : A^{p,q}(M) \to A^{p+q,r+s}(M)$ by
\[e(\omega) \eta = \omega \wedge \eta .\]
\end{Def}

\begin{Def}
We denote by $\mathscr{X}(M)$ the set of smooth vector fields on $M$. For $X \in \mathscr{X}(M)$ we define interior product operators by
\[
\begin{split}
i(X)& : A^{p,q}(M) \to A^{p-1,q}(M), \quad i(X) \omega = \omega (X,\dots;\dots),  \\
 \bar {i}(X)& : A^{p,q}(M) \to A^{p,q-1}(M), \quad \bar{i}(X) \omega = \omega ( \ \dots;X,\dots).
\end{split}
\]
\end{Def}

We denote by $v_g$ the volume form for $g$;
\[v_g= \sqrt { \det[g_{ij}]} \, dx^1 \wedge \dots \wedge dx^n.\]

\begin{Def}
We define $ \bigstar :A^{p,q}(M) \to A^{n-p,n-q}(M)$ by
\[ \omega \wedge \bigstar \eta = \langle \omega , \eta \rangle v_g \otimes \bar {v}_g , \quad \omega , \eta \in A^{p,q}(M),\]
where $ \langle \ , \ \rangle $ is a fiber metric on $ \wedge ^p T^* M \otimes \wedge ^q T^* M$ induced by $g$.
\end{Def}

Let $(E_1, \dots ,E_n)$ be a positive orthogonal frame field with respect to $g$ and$( \theta ^1, \dots , \theta ^n)$ be the dual frame field of  $(E_1, \dots ,E_n)$. Then we have
\[ \theta ^j = \bar {i}(E_j)g, \quad \bar { \theta }^j=i(E_j)g.\]
For a multi-index $I_p=(i_1, \dots ,i_p), \,  i_1< \dots <i_p$, we define $I_{n-p}=(i_{p+1}, \dots ,i_{n}), \, i_{p+1}< \dots <i_n$, where $(I_p , I_{n-p})$ is a permutation of $(1, \dots ,n)$. We denote by $ \epsilon_{(I_p,I_{n-p})}$ the signature of $(I_p , I_{n-p})$. Then by definition of $\bigstar$,
\[ \bigstar ( \theta ^{I_p} \otimes \overline { \theta ^{J_q}})= \epsilon _{(I_p ,I_{n-p})} \epsilon _{(J_q ,J_{n-q})} \theta ^{I_{n-p}} \otimes \overline { \theta ^{J_{n-q}}}. \]

\begin{Lem} \label{St}
\cite{Sh} The following identities hold on $A^{p,q}(M)$.  
\begin{enumerate}
\item[$(1)$] 
$\bigstar \bigstar =(-1)^{(p+q)(n+1)} $.
\item[$(2)$]
$ i(X)=(-1)^{p+1} \bigstar ^{-1}e( \bar {i}(X)g) \bigstar $,

$\bar {i}(X)=(-1)^{q+1} \bigstar ^{-1}e(i(X)g) \bigstar , \quad X \in \mathscr{X}(M)$.
\end{enumerate}
\end{Lem}

\begin{Lem} \label{IE}
\cite{Sh} The following equations hold for $ \omega \in A^{p,q}(M)$, $\eta \in A^{p-1,q}(M), \rho \in A^{p,q-1}(M)$ and $X \in \mathscr{X}(M)$.
\[ \langle i(X) \omega , \eta \rangle = \langle \omega ,e( \bar{i}(X)g) \eta \rangle ,\]
\[ \langle \bar{i}(X) \omega , \rho \rangle = \langle \omega ,e( i(X)g) \rho \rangle .\]
\end{Lem}

\begin{Def}
We define $L: A^{p,q}(M) \to A^{p+1,q+1}(M)$ and $\Lambda : A^{p,q}(M) \to A^{p-1,q-1}(M)$ by
\[\displaystyle L:=e(g)= \sum _j e( \theta ^j)e( \bar { \theta }^j),  \quad \Lambda := \sum _j i(E_j) \bar {i}(E_j).\]
\end{Def}

We obtain the following from Lemma \ref{IE}.

\begin{Cor} \label{L}
\cite{Sh} We have
\[ \langle \Lambda \omega, \eta \rangle = \langle \omega, L \eta \rangle, \quad \textrm{for} \ \omega \in A^{p,q}(M) \ \textrm{and} \ \eta \in A^{p-1,q-1}(M) .\]
\end{Cor}

\begin{Prop} \label{LL}
\cite{Sh} We have
\[[L, \Lambda ]=n-p-q, \quad \textrm{on} \ A^{p,q}(M).\]
\end{Prop}

\subsection{Differential operators for $(p,q)$-forms}

\begin{Def}
We define $\partial :A^{p,q}(M) \to A^{p+1,q}(M)$ and $\bar{\partial} :A^{p,q}(M) \to A^{p,q+1}(M)$ by
\[\partial =\sum_{i}e(dx^i)D_{\frac{\partial}{\partial x^i}}, \quad
\bar{\partial}=\sum_{i}e(\overline{dx^i})D_{\frac{\partial}{\partial x^i}}.\]
\end{Def}

Since $D$ is flat, we immediately obtain the following lemma.

\begin{Lem} \label{PP}
\cite{Sh} We have
\[\partial^2=0, \quad \bar{\partial}^2=0, \quad \partial \bar{\partial}=\bar{\partial}\partial.\]
\end{Lem}

We denote by $A^{p,q}(M,F)$ the space of $F$-valued $(p,q)$-forms. Since the transition functions of $\{s\}$ are constant, $\partial$ and $\bar{\partial}$ are extended on $A^{p,q}(M,F)$ by
\[\partial (s \otimes \omega) =s \otimes \partial \omega ,\]
\[\bar{\partial} (s \otimes \omega) =s \otimes \bar{\partial} \omega .\]

\begin{Def}
We denote by $A^{p,q}_0(M,F)$ the space of elements of $A^{p,q}(M,F)$ with compact supports. We define the inner product $(\ , \ )$on  $A^{p,q}_0(M,F)$ by
\[(\omega ,\eta)= \int_M\langle \omega , \eta \rangle _F  \ v_g,\]
where $\langle \ , \ \rangle$ is the metric on $F \otimes \wedge^p T^*M \otimes \wedge^q T^*M$ induced by $g$ and $h$. We set $ \| \omega \| = \sqrt {(\omega , \omega)}$.
\end{Def}

\begin{Def}
We define $A \in A^{1,0}(M)$ and $B \in A^{1,1}(M)$ by
\[A=- \partial \log h(s,s), \quad B= \bar{ \partial }A .\]
We call $A$ and $B$ the \emph{first Koszul form} and the \emph{second Koszul form} with respect to the fiber metric $h$, respectively.
\end{Def}

\begin{Rem}
Since the transition functions of $\{s\}$ are constant, $A$ and $B$ are globally well-defined.
\end{Rem}

\begin{Ex}
\cite{Sh} Let $\alpha$ and $\beta$ be the first Koszul form and the second Koszul form with respect to the Riemannian metric $g$, respectively. Then the the first Koszul form $A_K$ and the second Koszul form $B_K$ with respect to the fiber metric $g$ on $K$ are given by 
\[A_K=2\alpha , \quad B_K=2\beta.\]
\end{Ex}

\begin{Def} \label{Pr}
We define $\partial ^\prime _F :A^{p,q}(M,F) \to A^{p+1,q}(M,F)$ by
\[\partial ^\prime _F = \partial -e( A+ \alpha ).\]
We denote it by $\partial ^\prime $ if $(F,D^F ,h)$ is trivial.
\end{Def}

\begin{Thm} \label{Can}
We have
\[( \partial ^\prime _F)^2=0, \quad \partial ^\prime _F \bar{\partial}- \bar{\partial} \partial ^\prime _F=e(B+ \beta) .\]
\end{Thm}

\begin{Prf}
We obtain
\[
\begin{split}
( \partial ^\prime _F)^2
&=(\partial -e( A+ \alpha ))(\partial -e( A+ \alpha )) \\
&=\partial ^2-e(\partial (A+ \alpha))+e(A+ \alpha) \partial -e(A+ \alpha) \partial +e( A+ \alpha )e( A+ \alpha ) \\
&=0,
\end{split}
\]
and
\[
\begin{split}
\partial ^\prime _F \bar{\partial}
&=(\partial -e(A+ \alpha)) \bar{\partial}
= \partial \bar{\partial} -e(A+ \alpha) \bar{\partial},\\
\bar{\partial} \partial ^\prime _F
&= \bar{\partial} (\partial -e(A+ \alpha))
= \bar{\partial} \partial -e (B+ \beta))-e(A+ \alpha) \bar{\partial}.
\end{split}
\]
Hence
\[\partial ^\prime _F \bar{\partial}- \bar{\partial} \partial ^\prime _F=e(B+ \beta).\]
\begin{flushright}
$\Box$
\end{flushright}
\end{Prf}

\begin{Def} \label{di}
We define $\delta ^\prime _F$, $\delta_F : A^{p,q}(M,F) \to A^{p-1,q}(M,F)$ and $\bar{\delta}_F : A^{p,q}(M,F) \to A^{p,q-1}(M,F)$ by
\[\delta ^\prime _F = (-1)^p \bigstar ^{-1} \partial \bigstar , \quad 
\delta _F = (-1)^p \tilde{\bigstar}_F ^{-1} \partial \tilde{\bigstar} _F , \quad
\bar{\delta} _F=(-1)^q \tilde{\bigstar}_F ^{-1} \bar{\partial} \tilde{\bigstar} _F .\] 
We denote them by $ \delta ^\prime$, $\delta$ and $\bar{\delta}$ if $(F,D^F ,h)$ is trivial.
\end{Def}

\begin{Prop} \label{Ad}
\cite{Sh} The operators $\delta _F$ and $\bar{\delta} _F$ are the adjoint operators of  $\partial$ and $\bar{\partial}$ with respect to the inner product $( \ , \ )$ respectively, that is, for $ \omega \in A^{p,q}(M,F)$, $\eta \in A^{p-1,q}_0(M,F)$ and $\rho \in A^{p,q-1}_0(M,F)$ we have
\[ ( \delta _F \omega , \eta)=( \omega , \partial \eta), \ ( \bar{\delta} _F \omega , \eta)=( \omega , \bar{\partial} \eta).\]
\end{Prop}

\begin{Cor} \label{Adp}
We have
\[ \delta _F ^\prime = \delta _F - i(X_{A+ \alpha}),\]
where $ \bar {i}(X_{A+ \alpha})g=A+ \alpha$. In addition, $ \delta ^\prime _F$ is the adjoint operator of $ \partial ^\prime _F$ with respect to the inner product $( \ , \ )$.
\end{Cor}

\begin{Prf}
On $A^{p,q}(M,F)$ we have
\[
\begin{split}
\delta _F &=(-1)^p \tilde{\bigstar} _F ^{-1} \partial \tilde{\bigstar} _F  \\
&=(-1)^p h(s,s)^{-1} ( \det[g_{ij}])^{\frac{1}{2}} \bigstar ^{-1} \partial (h(s,s)( \det[g_ij])^{- \frac{1}{2}} \bigstar) \\
&=(-1)^p \bigstar ^{-1} \partial \bigstar +(-1)^{p+1} \bigstar ^{-1} e(- \partial \log (h(s,s) (\det[g_{ij}])^{- \frac{1}{2}})) \bigstar \\
&= \delta _F ^\prime +(-1)^{p+1} \bigstar ^{-1} e(A+ \alpha) \bigstar .
\end{split}
\]
Hence, it follows  from Lemma \ref{St} that we have
\[ \delta _F ^\prime = \delta _F - i(X_{A+ \alpha}),\]
where $ \bar {i}(X_{A+ \alpha})g=A+ \alpha$.
Therefore, $\delta_F^\prime$ is the adjoint operator of $ \partial _F ^\prime = \partial - e(A+ \alpha)$ from Proposition \ref{Ad} and Lemma \ref{IE},.
\begin{flushright}
$\Box$
\end{flushright}
\end{Prf}

\begin{Def}
We define the connection $ \mathscr{D}$ and $ \overline {\mathscr{D}}$ on $ \wedge ^p T^*M \otimes \wedge ^q T^*M$ as follows:
For $\omega \in A^p(M) $ and $\eta \in A^q(M),X \in \mathscr{X}(M)$
\[ 
\begin{split}
\mathscr{D} _X (\omega \otimes \bar{\eta}) &=2 \gamma _X \omega \otimes \bar{\eta} + D_X ( \omega \otimes \bar{\eta} ), \\
\overline{\mathscr{D}}_X (\omega \otimes \bar{\eta}) &=2 \omega \otimes \overline { \gamma _X \eta} + D_X ( \omega \otimes \bar{\eta} ) ,
\end{split}
\]
where $ \gamma = \nabla  - D $ and $ \nabla $ is the Levi-Civita connection of $g$ （Definition \ref{na}）．
\end{Def}

The following lemma follows from Proposition \ref{H}.

\begin{Lem} \label{He}
\cite{Sh} The following conditions are equivalent.
\item[$(1)$] $(D,g)$ is a Hessian structure．
\item[$(2)$] $ \partial g =0 \quad ( \Leftrightarrow \bar{\partial} g=0)$．
\item[$(3)$] $ \mathscr{D} g=0 \quad (\Leftrightarrow \overline {\mathscr{D}} g =0)$．
\end{Lem}

Let $D^*$ be the dual connection of $D$ with respect to $g$ （Definition \ref{Du}）. We obtain the following from Proposition \ref{H} and \ref{S}

\begin{Lem} \label{scD}
Let $(D,g)$ be a Hessian structure. Then we have
\[
\begin{split}
\mathscr{D}_X (\omega \otimes \bar{\eta})&=D^*_X \omega \otimes \bar{\eta} + \omega \otimes \overline{D_X \eta} ,\\
\overline{\mathscr{D}}_X (\omega \otimes \bar{\eta})&=D_X \omega \otimes \bar{\eta} + \omega \otimes \overline{D^*_X \eta} ,
\end{split}
\]
for $\omega \in A^p(M)$ and $\eta \in A^{q}(M),X \in \mathscr{X}(M)$.
\end{Lem}

\begin{Prop} \label{PD}
\cite{Sh} Let $(D,g)$ is a Hessian structure. Then we have
\[ \partial = \sum _j e( \theta ^j ) \mathscr{D}_{E_j}, \quad  
\bar{\partial} = \sum _j e( \bar{\theta} ^j ) \overline{\mathscr{D}}_{E_j} .\]
\end{Prop}

\begin{Prop} \label{del}
\cite{Sh} Let $(D,g)$ is a Hessian structure. Then we have
\[
\begin{split}
\delta ^\prime _F &= - \sum _j i(E_j) \overline{\mathscr{D}}_{E_j}, \\
\bar{\delta} _F &= - \sum _j \bar{i}(E_j) \mathscr{D}_{E_j} + \bar{i}(X_{A+ \alpha}),
\end{split}
\]
where $\bar{i}(X_{A+\alpha})g = A+\alpha$.
\end{Prop}

The following theorem is an analogue of K\"{a}hler identities.

\begin{Thm} \label{KI}
Let $(D,g)$ is a Hessian structure. Then we have
\[ \Lambda \partial ^\prime _F + \partial ^\prime _F \Lambda = - \bar{\delta} _F , \quad 
\Lambda \bar{\partial} + \bar{\partial} \Lambda = -\delta ^\prime _F ,\]
\[L \delta ^\prime _F + \delta ^\prime _F L = - \bar{\partial} , \quad 
L\bar{\delta} _F + \bar{\delta} _F L = - \partial ^\prime _F .\]
\end{Thm}

\begin{Prf}
It follows from Proposition \ref{del}, \ref{He} and \ref{PD} that we obtain
\[
\begin{split}
\delta ^\prime _F L
&= -\sum _j i(E_j) \overline{\mathscr{D}}_{E_j} L
= -\sum _j i(E_j)L \overline{\mathscr{D}}_{E_j} \\
&=- \sum _{j,k} i(E_j) e(\theta^k)e(\bar{\theta}^k) \overline{\mathscr{D}}_{E_j} \\
&=- \sum _{j,k} e(\bar{\theta}^k)(\delta_j^k - e(\theta^k)i(E_j)) \overline{\mathscr{D}}_{E_j} \\
&=- \sum _j e(\bar{\theta}^j) \overline{\mathscr{D}}_{E_j} + \sum _k e(\bar{\theta}^k)e(\theta^k) \sum _j i(E_j) \overline{\mathscr{D}}_{E_j} \\
&=- \bar{\partial} - L \delta ^\prime _F.
\end{split}
\]
Similarly, we have
\[- \sum _j \bar{i}(E_j) \mathscr{D}_{E_j} L =- \partial + L \sum_j \bar{i}(E_j) \mathscr{D}_{E_j}.\]
Moreover,
\[
\begin{split}
\bar{i}(X_{A+\alpha})L
&=\bar{i}(X_{A+\alpha}) \sum _k e(\theta^k)e(\bar{\theta}^k) \\
&=\sum _k e(\theta^k) \{(A+\alpha)(E_k) - e(\bar{\theta}^k) \bar{i}(X_{A+\alpha}) \} \\
&=e(A+\alpha) - L \bar{i}(X_{A+\alpha}).
\end{split}
\]
Hence it follows from Corollary \ref{Ad} that
\[\bar{\delta}_F L
=(- \sum_j \bar{i}(E_j) \mathscr{D}_{E_j} + \bar{i}(X_{A+\alpha})) L
=-\partial ^\prime _F - L \bar{\delta} _F.\]
We have the other equalities by taking the adjoint operators.
\begin{flushright}
$\Box$
\end{flushright}
\end{Prf}

\begin{Def}
We define the Laplacians $ \square ^\prime _F$ and $\bar { \square }$ with respect to $ \partial ^\prime _F$ and $ \bar { \partial }$ by
\[ \square ^\prime _F = \partial ^\prime _F \delta ^\prime _F +  \delta ^\prime _F \partial ^\prime _F , \quad \bar { \square }_F = \bar { \partial } \delta _F + \delta _F \bar { \partial } .\]
We denote them by $\square ^ \prime$ and $\bar { \square }$ if $(F,D^F,h)$ is trivial.
\end{Def}

The following theorem is an analogue of Kodaira-Nakano identity.

\begin{Thm} \label{Id}
Let $(D,g)$ is a Hessian structure. Then we have
\[ \bar { \square }_F= \square ^\prime _F +[ e( \beta +B ), \Lambda ]. \]
\end{Thm}

\begin{Prf}
It follows from Theorem \ref{Can} and \ref{KI}  that we obtain
\[
\begin{split}
\bar{\square}_F &= \bar{\partial} \bar{\delta}_F + \bar{\delta}_F \bar{\partial}
=-\bar{\partial} (\Lambda \partial ^\prime _F + \partial ^\prime _F \Lambda)  - (\Lambda \partial ^\prime _F + \partial ^\prime _F \Lambda) \bar{\partial} \\
&=(\Lambda \bar{\partial} + \delta ^\prime _F) \partial ^\prime _F - \bar{\partial}  \partial ^\prime _F \Lambda - \Lambda \partial ^\prime _F \bar{\partial} + \partial ^\prime _F (\bar{\partial} \Lambda +\delta ^\prime _F) \\
&=\delta ^\prime _F \partial ^\prime _F + \partial ^\prime _F \delta ^\prime _F +(\partial ^\prime _F \bar{\partial} - \bar{\partial} \partial ^\prime _F) \Lambda - \Lambda (\partial ^\prime _F \bar{\partial} - \bar{\partial} \partial ^\prime _F )  \\
&=\square _F ^\prime + [e(B+\beta),\Lambda].
\end{split}
\]
\begin{flushright}
$\Box$
\end{flushright}
\end{Prf}

The following theorem is an analogue of Kodaira-Nakano vanishing theorem.

\begin{Thm} \label{C}
\cite{Sh} Let $(M,D)$ be an oriented $n$-dimensional compact flat manifold and $(F,D^F)$ be a flat line bundle over $M$.  We set
\[H^{p,q}_{\bar{\partial}}(M,F)= \frac { \textrm{Ker} \ [ \ \bar { \partial }:A^{p,q}(M,F) \to A^{p,q+1}(M,F) \ ]} { \textrm{Im} \ [ \ \bar { \partial }:A^{p,q-1}(M,F) \to A^{p,q}(M,F) \ ]}.\]
Assume there exists a fiber metric $h$ on $F$ such that $B+\beta > 0$, where $B$ and $\beta$ are the second Koszul form with respect to $h$ and $g$ respectively. Then we have
\[H^{p,q}_{\bar{\partial}}(M,F)=0, \quad \textrm{for} \ p+q>n.\]
\end{Thm}

\section{Vanishing theorems of $L^2$-cohomology groups} \label{L2}

\subsection{$L^2$-cohomology groups}

\begin{Def}
We denote $L^{p,q}(M,g,F,h)$ by the completion of $A^{p,q}_0(M,F)$ with respect to the $L^2$-inner product $(\ , \ )$. The space $L^{p,q}(M,g,F,h)$ is identified with the space of square-integrable sections of $F \otimes \wedge^p T^*M \otimes \wedge^q T^*M$.
\end{Def}

\begin{Def} \label{Ge}
For $\omega \in L^{p,q}(M,g,F,h)$ we define $\bar{\partial} \omega$ and $\bar{\delta}_F \omega$ as follows:
\[(\bar{\partial} \omega , \eta)=(\omega , \bar{\delta}_F \eta), \quad \textrm{for} \ \eta \in A^{p,q+1}_0(M,F),\]
\[(\bar{\delta}_F \omega , \rho )=(\omega , \bar{\partial} \rho), \quad \textrm{for} \ \rho \in A^{p,q-1}_0(M,F).\]
In general, we cannot say $\bar{\partial} \omega \in L^{p,q+1}(M,F)$ and $\bar{\delta}_F \omega \in L^{p,q-1}(M,g,F,h)$. We set
\[
\begin{split}
W^{p,q}(M,g,F,h)&= \{\omega \in L^{p,q}(M,g,F,h) \mid \bar{\partial} \omega \in L^{p,q+1}(M,g,F,h), \bar{\delta}_F \omega \in L^{p,q-1}(M,g,F,h) \}, \\
D^{p,q}(M,g,F,h)&=\{\omega \in L^{p,q}(M,g,F,h) \mid \bar{\partial} \omega \in L^{p,q+1}(M,g,F,h) \}.
\end{split}
\]
In addition, we define the norm $\| \ \|_W$ on $W^{p,q}(M,g,F,h)$ by
\[ \| \omega \|_W = \| \omega \| + \| \bar{\partial} \omega \| + \| \bar{\delta}_F \omega \|, \quad \omega \in W^{p,q}(M,g,F,h) .\]
The space $W^{p,q}(M,g,F,h)$ is complete with respect to $\| \ \|_W$.
\end{Def}

\begin{Prop} \label{De}
\cite{Je} If $g$ is complete, the space $A^{p,q}_0(M,F)$ is dense in $W^{p,q}(M,g,F,h)$ with respect to $\| \ \|_W$.
\end{Prop}

\begin{Def}
We define the  \emph{$L^2$-cohomology group} of $(p,q)$-type by
\[L^2H^{p,q}_{\bar{\partial}}(M,g,F,h)= \frac { \textrm{Ker} \ [ \ \bar { \partial }:D^{p,q}(M,g,F,h) \to D^{p,q+1}(M,g,F,h) \ ]} {\overline {\textrm{Im} \ [ \ \bar { \partial }:D^{p,q-1}(M,g,F,h) \to D^{p,q}(M,g,F,h) \ ]}},\]
where $\overline {\textrm{Im} \ [ \ \bar { \partial }:D^{p,q-1}(M,g,F,h) \to D^{p,q}(M,g,F,h) \ ]}$ is the closure of $\textrm{Im} \ [ \ \bar { \partial }:D^{p,q-1}(M,g,F,h) \to D^{p,q}(M,g,F,h) \ ]$ with respect to the $L^2$-norm $\| \ \|$.
\end{Def}

\subsection{Vanishing theorems of Kodaira-Nakano type }

\begin{Lem} \label{b}
Assume $B+\beta$ is positive definite. For the eigenvalues $\lambda_1 \le \dots \le \lambda_n$ of the matrix $\displaystyle \Bigl[ \sum _k g^{ik}( \beta +B)_{kj} \Bigr]$, we set $\displaystyle b_q=\sum_{j=1}^q \lambda_j$. Then  we have
\[\| \bar{\partial} \omega \|^2 + \| \bar{\delta}_F \omega\|^2 \ge \| b_q^{\frac{1}{2}} \omega \|^2 , \quad \textrm{for} \ \omega \in A^{n,q}_0(M,F).\]
\end{Lem}

\begin{Prf}
By Theorem \ref{Id} we obtain
\[\| \bar{\partial} \omega \|^2 + \| \bar{\delta}_F \omega\|^2
=(\bar{\square}_F \omega , \omega) 
= (\square^\prime_F \omega , \omega) + ([e(B+\beta),\Lambda] \omega , \omega)
\ge  ([e(B+\beta),\Lambda] \omega , \omega) .\]
Hence it is sufficient to show$ ([e(B+\beta),\Lambda] \omega , \omega) \ge \| b_q^{\frac{1}{2}} \omega \|^2$.

We take the positive orthonormal frame field $\{E_1, \dots ,E_n\}$ on $TM$ where the matrix $[(B+\beta)(E_i,E_j)]$ is diagonal. We set $\mu_j =(B+\beta)(E_j,E_j)$. Using the dual frame field $\{ \theta^1, \dots ,\theta^n \}$ of $\{E_1, \dots ,E_n\}$ , we denote $\omega \in A^{n,q}_0(M,F)$ by
\[\omega = \sum_{J_q} \omega_{J_q} \otimes \bar{\theta}^{J_q} , \quad  \ \omega_{J_q} \in A^n_0(M,F).\]
Hence
\[
\begin{split}
[e(B+\beta),\Lambda]\omega &= e(B+\beta) \Lambda \omega \\
&=\sum_j \mu_j e(\theta^j)e(\bar{\theta}^j) \sum _k i(E_k)\bar{i}(E_k) \sum_{J_q} \omega_{j_q} \otimes \bar{\theta}^{J_q} \\
&=\sum_{j,J_q} \mu_j \omega_{J_q} \otimes e(\bar{\theta}^j) \bar{i}(E_j)\bar{\theta}^{J_q} \\
&=\sum_{J_q} \sum_{j \in J_q} \mu_j \omega_{J_q} \otimes \bar{\theta}^{J_q}.
\end{split}
\]
Therefore
\[
\begin{split}
([e(B+\beta),\Lambda]\omega , \omega)
&=\int_M \sum_{J_q} \sum_{j \in J_q} \mu_j \langle \omega_{J_q} , \omega_{J_q} \rangle_h v_g \\
&\ge \int_M \sum_{J_q} b_q \langle \omega_{J_q} , \omega_{J_q} \rangle_h v_g 
=\| b_q^{\frac{1}{2}} \omega \|^2.
\end{split}
\]
\begin{flushright}
$\Box$
\end{flushright}
\end{Prf}

\begin{Main} \label{L}
Let $(M,D,g)$ be an oriented $n$-dimensional complete Hessian manifold and $(F,D^F)$ a flat line bundle over $M$. We denote by $h$ a fiber metric on $F$. Assume  $B+ \beta$ is positive definite, where $B$ and $\beta$ are the second Koszul forms with respect to fiber metric $h$ and Hessian metric $g$ respectively. For $q \ge 1$ let $b_q$ be the same as in Lemma \ref{b}. Then for all $v \in L ^{n,q}(M,g,F,h)$ such that $ \bar { \partial }v=0$ and  $b_q^{-\frac{1}{2}}v \in L^{p,q}(M,g,F,h)$, there exists $u \in L^{n,q-1}(M,g,F,h)$ such that 
\[ \bar { \partial }u=v, \quad \| u \| \le \| b_q^{- \frac {1}{2}} v \| .\]
In particular, if there exists $\varepsilon > 0$ such that $B + \beta - \varepsilon g$ is positive definite, we have 
\[L^2H^{n,q}_{\bar{\partial}}(M,g,F,h)=0, \quad \textrm{for} \ q \geq 1.\] 
\end{Main}

\begin{Prf}
We set $\textrm{Ker} \bar{\partial}=\{\omega \in L^{n,q}(M,g,F,h) \mid \bar{\partial} \omega =0\}$. Since $\textrm{Ker} \bar{\partial}$ is a closed subspace in $ L^{n,q}(M,g,F,h)$, we have 
\[ L^{n,q}(M,g,F,h) = \textrm{Ker}\bar{\partial} \oplus (\textrm{Ker}\bar{\partial})^\perp, \]
where  $(\textrm{Ker}\bar{\partial})^\perp$ is the orthogonal complement of $\textrm{Ker} \bar{\partial}$. We denote $\omega \in L^{n,q}(M,g,F,h)$ by
\[\omega=\omega_1 + \omega_2, \quad \omega_1 \in \textrm{Ker}\bar{\partial}, \quad \omega_2 \in (\textrm{Ker}\bar{\partial})^\perp .\]
For $\eta \in A^{n,q-1}_0(M,F)$, we have
\[(\bar{\delta}_F \omega , \eta)=(\omega , \bar{\partial}\eta)=0,\]
and so
\[\bar{\delta}_F \omega_2=0.\]

Since $v \in \textrm{Ker} \bar{\partial}$ by assumption, we obtain
\[|(v,\omega)|^2=|(v,\omega_1)|^2=|(b_q^{-\frac{1}{2}}v,b_q^{\frac{1}{2}} \omega_1)|^2
\le \|b_q^{-\frac{1}{2}}v\|^2 \|b_q^{\frac{1}{2}} \omega_1\|^2.\]

Assume $\omega \in W^{n,q}(M,g,F,h)$. Then
\[\bar{\partial}\omega_1=0, \quad \bar{\delta}_F \omega_1 = \bar{\delta}_F \omega \in L^{p,q}(M,g,F,h), \]
and so $\omega_1 \in W^{n,q}(M,g,F,h)$. Hence by Proposition \ref{De}, $\omega_1$ satisfies the inequality in Lemma \ref{b}:
\[\|b_q^{\frac{1}{2}} \omega_1\|^2 
\le \| \bar{\partial} \omega_1 \|^2 + \| \bar{\delta}_F \omega_1\|^2
=\| \bar{\delta}_F \omega_1\|^2 
=\| \bar{\delta}_F \omega \|^2 < \infty  .\]
Therefore for $\omega \in W^{n,q}(M,g,F,h)$ we have
\[|(v,\omega)|^2 \le \|b_q^{-\frac{1}{2}}v\|^2  \| \bar{\delta}_F \omega \|^2 < \infty. \]
By this inequality a linear functional
$\lambda : \bar{\delta}_F W^{n,q}(M,g,F,h) \ni \bar{\delta}_F \omega \mapsto (v, \omega) \in \mathbb{R}$ is well-defined and the operator norm $C$ is
\[C \le \|b_q^{-\frac{1}{2}}v\| < \infty.\]

We set $\textrm{Ker} \bar{\delta}_F=\{\omega \in L^{n,q}(M,g,F,h) \mid \bar{\delta}_F \omega =0\}$. $\textrm{Ker} \bar{\delta}_F$ is also a closed subspace in $ L^{n,q}(M,g,F,h)$ and 
\[ L^{n,q}(M,g,F,h) = \textrm{Ker}\bar{\delta}_F \oplus (\textrm{Ker}\bar{\delta}_F)^\perp ,\]
where $(\textrm{Ker}\bar{\delta}_F)^\perp$ is the orthogonal complement of $\textrm{Ker}\bar{\delta}_F$. In the same way we have $(\textrm{Ker}\bar{\delta}_F)^\perp \subset \textrm{Ker}\bar{\partial}$ and for $\hat{\omega} \in  (\textrm{Ker}\bar{\delta}_F)^\perp \cap W^{n,q}(M,g,F,h)$,
\[\|b_q^{\frac{1}{2}} \hat{\omega} \|^2 
\le \| \bar{\partial} \hat{\omega} \|^2 + \| \bar{\delta}_F\hat{ \omega} \|^2
=\| \bar{\delta}_F \hat{\omega} \|^2 .\]
Let $\{\eta_k\} \subset \bar{\delta}_F W^{n,q}(M,g,F,h)$ be a Cauchy sequence with respect to the norm $\| \ \|$ on $L^{n,q-1}(M,g,F,h)$. Each $\eta_k$ is denoted by
\[ \eta_k = \bar{\delta}_F \hat{\omega}_k , \quad \hat{\omega}_k \in (\textrm{Ker}\bar{\delta}_F)^\perp \cap W^{n,q}(M,g,F,h) ,\]
and by the said inequality $\{ \hat{\omega}_k \}$ is also a Cauchy sequence with respect to the norm $\| \ \|$ on $L^{n,q}(M,g,F,h)$. This implies $\{ \hat{\omega}_k \}$ is a Cauchy sequence with respect to the norm $\| \ \|_W$ on $W^{n,q}(M,g,F,h)$. Hence by completeness of $W^{n,q}(M,g,F,h)$ with respect to $\| \ \|_W$, we have \[\hat{\omega}_k \to \hat{\omega} \in W^{n,q}(M,g,F,h) \quad (k \to \infty),\]
and 
\[\eta_k \to \bar{\delta}_F \hat{\omega} \quad (k \to \infty).\]
Therefore, $\bar{\delta}_F W^{n,q}(M,g,F,h)$ is a closed space of $L^{n,q-1}(M,g,F,h)$ with respect to the norm $\| \ \|$.

From the above, by applying Riesz representation theorem to the linear functional $\lambda: \bar{\delta}_F W^{n,q}(M,g,F,h) \to \mathbb{R}$, there exists $u \in \bar{\delta}_F W^{n,q}(M,g,F,h)$ such that
\[
\begin{cases}
\lambda(\eta) = (u,\eta), \quad \eta \in \bar{\delta}_F W^{n,q}(M,g,F,h) \\ 
\|u\| = C  \le \|b_q^{-\frac{1}{2}}v\|.
\end{cases}
\]
By the first equation, for all $\omega \in A^{n,q}_0(M,F)$ we have 
\[(v,\omega)=\lambda(\bar{\delta}_F \omega)=(u, \bar{\delta}_F \omega),\] 
and so
\[\bar{\partial} u = v.\]
This implies the first assertion.

Suppose there exists $\varepsilon > 0$ such that $B + \beta - \varepsilon g$ is positive definite. Then by the definition of $b_q$
$b_q \ge \varepsilon q$. Hence for all $v \in L^{n,q}(M,g,F,h)$ we obtain
\[ \int_M \langle b_q^{-\frac{1}{2}},b_q^{-\frac{1}{2}} \rangle v_g \le (\varepsilon q)^{-1}\int_M \langle v,v \rangle v_g < \infty,\]
that is,
\[b_q^{-\frac{1}{2}}v \in L^{n,q}(M,g,F,h).\]
This implies the second assertion.
\begin{flushright}
$\Box$
\end{flushright}
\end{Prf}

The following theorem is an extension of Theorem \ref{C}

\begin{Main} \label{LC}
Let $(M,D,g)$ be an oriented $n$-dimensional complete Hessian manifold and $(F,D^F)$ a flat line bundle over $M$. We denote by $h$ a fiber metric on $F$. Assume that there exists $\varepsilon >0$ such that $B+ \beta= \varepsilon g$ where $B$ and $\beta$ are the second Koszul forms with respect to fiber metric $h$ and Hessian metric $g$ respectively. Then if for $p+q>n$ and all $v \in L ^{p,q}(M,g,F,h)$ such that $ \bar { \partial }v=0$, there exists $u \in L^{p,q-1}(M,g,F,h)$ such that 
\[ \bar { \partial }u=v, \quad  \| u \| \le  \{ \varepsilon(p+q-n)\}^{- \frac {1}{2}} \| v \| .\] 
In particular, we have
 \[L^2H^{p,q}_{\bar{\partial}}(M,g,F,h)=0, \quad \textrm{for} \ p+q>n.\]
\end{Main}

\begin{Prf}
By Proposition \ref{LL}, on $A^{p,q}(M,F)$ we have 
\[[e(B+\beta), \Lambda]=\varepsilon[L,\Lambda]=\varepsilon (p+q-n).\]  Hence by Theorem \ref{Id}, for all $\omega \in A^{p,q}_0(M,F)$ we obtain
\[\| \bar{\partial} \omega \|^2 + \| \bar{\delta}_F \omega\|^2 \ge \varepsilon (p+q-n) \| \omega \|^2.\]
Then the assertions are proved similarly with Main theorem \ref{L}.
\begin{flushright}
$\Box$
\end{flushright}
\end{Prf}

\begin{Cor} \label{Rn}
Let $(\mathbb{R}^n,g)$ be a Euclidean space，$D$ be a canonical affine connection on  $\mathbb{R}^n$, and $(F=\mathbb{R}^n  \times \mathbb{R},D^F)$ be a trivial flat line bundle on $\mathbb{R}^n$. In addition, we define a fiber metric $h$ on $F$ by 
\[ h(s,s)=e^{-\varphi},\] 
where $\displaystyle \varphi(x)=\frac{1}{2} \sum_i (x^i)^2$ and $s: \mathbb{R}^n \ni x \mapsto (x,1) \in F$.
Then for $q \ge 1$ and $v \in L ^{p,q}(\mathbb{R}^n,g,F,h)$ such that $\bar { \partial }v=0$, there exists $u \in L^{p,q-1}(\mathbb{R}^n,g,F,h)$ such that 
\[ \bar { \partial }u=v, \quad \| u \| \le q^{- \frac {1}{2}} \| v \|.\]
In particular, we have 
\[L^2H^{p,q}_{\bar{\partial}}(\mathbb{R}^n,g,F,h)=0, \quad \textrm{for} \ p \geq 0 \ \textrm{and} \ q \geq 1.\]
\end{Cor}

\begin{Prf}
The Hessian metric $g=Dd\varphi$ is complete and the second Koszul forms with respect to $h$ and $g$ are 
\[B=-\partial \bar{\partial} \log h(s,s)=\partial \bar{\partial} \varphi=g, \quad  \beta=\frac{1}{2} \partial \bar{\partial}\det [\delta_{ij}]=0.\]
Hence by Main theorem \ref{LC}, for $p=n$ we obtain the assertion. 

Next, we consider the case of $p=0$. For $v \in L^{0,q}(\mathbb{R}^n,g,F,h)$ we set
\[\hat{v}=dx^1 \wedge \dots \wedge dx^n \otimes v\]
Then we have $\hat{v} \in L^{n,q}(\mathbb{R}^n,g,F,h)$ and $\| \hat{v} \|=\|v\|$. Since $\bar{\partial} v=0$ and $\bar{\partial} \hat{v}=0$ are equivalent, by Main theorem \ref{LC} there exists $\hat{u} \in L^{n,q-1}(\mathbb{R}^n,g,F,h)$ such that $ \bar { \partial }\hat{u}=\hat{v}$ and $ \| \hat{u} \| \le q^{- \frac {1}{2}} \| \hat{v}\| $. Here $\hat{u}$ is denoted by
\[\hat{u}=dx^1 \wedge \dots \wedge dx^n \otimes u, \quad u \in L^{0,q-1}(\mathbb{R}^n,g,F,h),\]
and so
\[dx^1 \wedge \dots \wedge dx^n \otimes \bar{\partial} u
=\bar{\partial} \hat{u} =\hat{v} =dx^1 \wedge \dots \wedge dx^n \otimes v .\]
Therefore, we have $\bar{\partial}u=v.$ Moreover, we obtain
\[\| u \| = \| \hat{u} \| \le q^{-\frac{1}{2}}\| \hat{v} \| = q^{-\frac{1}{2}} \| v \|.\]
Hence the assertion for $p=0$ follows.

Finally, for $p \ge 1$，$v \in L^{p,q}(\mathbb{R}^n,g,F,h)$ is denoted by
\[v=\sum_{I_p}dx^{I_p} \otimes v_{I_p}, \quad I_p=(i_1, \dots ,i_p), \quad 1 \le i_1 < \dots < i_p \le n, \quad v_{I_p} \in L^{0,q}(\mathbb{R}^n,g,F,h), \]
and we have
\[\| v \|^2 = \sum_{I_p} \| v_{I_p} \|^2.\]
If $\bar{\partial}v=0$, for all $I_p$ we obtain $\bar{\partial}v_{I_p}=0$. Hence by the case of $p=0$, there exists $\{u_{I_p}\} \subset L^{0,q-1}(\mathbb{R}^n,g,F,h)$ such that $\bar { \partial} u_{I_p}=v_{I_p}$ and $ \| u_{I_p} \| \le q^{- \frac {1}{2}} \| v_{I_p} \| $. Here we set
\[u=\sum_{I_p}dx^{I_p} \otimes u_{I_p}.\]
Then we have
\[\bar{\partial}u=\sum_{I_p}dx^{I_p} \otimes \bar{\partial}u_{I_p}
=\sum_{I_p}dx^{I_p} \otimes v_{I_p}=v,\]
\[\|u\|^2=\sum_{I_p} \|u_{I_p}\|^2 \le \sum_{I_p} q^{-1} \|v_{I_p}\|^2 = q^{-1} \|v\|^2.\]
This completes the proof. 
\begin{flushright}
$\Box$
\end{flushright}
\end{Prf}

\begin{Cor} \label{Reg}
Let $\Omega \in \mathbb{R}^n$ be a regular convex domain, $D$ be a canonical affine connection on $\Omega$, $g$ be a Hessian metric defined by Theorem \ref{M}. Then for  $p+q>n$ and $v \in L ^{p,q}(\Omega,g)$ such that $\bar { \partial }v=0$, there exists $u \in L^{p,q-1}(\Omega,g)$ such that 
\[ \bar { \partial }u=v, \quad  \| u \| \le (p+q-n)^{- \frac {1}{2}} \| v \| .\]
In particular, we have 
\[L^2H^{p,q}_{\bar{\partial}}(\Omega,g)=0, \quad \textrm{for} \ p+q>n.\]
\end{Cor}

\begin{Prf}
Since $g$ is complete and $\beta=g$, the assertion follows from Main theorem \ref{LC}.
\begin{flushright}
$\Box$
\end{flushright}
\end{Prf}

Let $\Omega \in \mathbb{R}^{n-1}$ be a regular convex domain and we set
$V=\{(ty,t) \in \mathbb{R}^n \mid y \in \Omega,t>0\}$. Let $\tilde{D}$ be a canonical affine connection on $V$ and $\tilde{g}$ be a Hessian metric on $(V,\tilde{D})$ defined by Theorem \ref{M}. In addition, we define an action $\rho : \mathbb{Z} \to \textrm{GL}(V)$ by
\[\rho(k)x=e^k x, \quad k \in \mathbb{Z}, \quad x \in V .\]
Then we have $\mathbb{Z} \backslash V \simeq \Omega \times S^1$. Moreover, this action preserves $(\tilde{D},\tilde{g})$ and so a Hessian structure $(D,g)$ on $\Omega \times S^1$ is defined by projecting $(\tilde{D},\tilde{g})$ on $\Omega \times S^1$. The Hessian metric $g$ is complete and the second Koszul form with respect to $g$ is equal to $g$. Hence the following theorem follows from Main theorem \ref{LC}.

\begin{Cor} \label{ReS}
Let $(\Omega \times S^1,D,g)$ be as above. Then for $p+q>n$ and  $v \in L ^{p,q}(\Omega \times S^1,g)$ such that $\bar { \partial }v=0$, there exists $u \in L^{p,q-1}(\Omega \times S^1,g)$ such that
\[ \bar { \partial }u=v, \quad \| u \| \le (p+q-n)^{- \frac {1}{2}} \| v \|.\]
In particular, we have 
\[L^2H^{p,q}_{\bar{\partial}}(\Omega \times S^1,g)=0, \quad \textrm{for} \ p+q>n.\]
\end{Cor}

\subsection{$L^2$-cohomology groups on regular convex cones}

\begin{Def}
A regular convex domain $\Omega$ in $\mathbb{R}^n$ is said to be a \emph{regular convex cone} if, for any $x$ in $\Omega$ and any positive real number $\lambda$, $\lambda x$ belongs to $\Omega$.
\end{Def}

\begin{Thm}
Let $(\Omega,D,g=Dd\varphi)$ be a regular convex cone in $\mathbb{R}^n$ with the Cheng-Yau metric (Theorem \ref{M}). Then we have the following equations.
\item[$(1)$] $\displaystyle \sum_j x^j \frac{\partial \varphi}{\partial x^j} = -n$.
\item[$(2)$] $\displaystyle\textrm{grad} \, \varphi = -\sum_j x^j \frac{\partial}{\partial x^j}$.
\item[$(3)$] $\displaystyle \sum_k x^k\gamma_{ijk} = -g_{ij}$.
\end{Thm}

\begin{Prf}
By the proof of Corollary \ref{A}, for $t>0$ and $x \in \Omega$ we have
\[\varphi (tx) = \varphi (x) - n \log t  .\]
Then we obtain
\[\sum_j x^j \frac{\partial \varphi}{\partial x^j} 
= \left. \frac{d}{dt} \right|_{t=1} \varphi (tx)
= -n.\]
Taking the derivative of both sides with respect to $x^i$ we have
\[(*) \quad \frac{\partial \varphi}{\partial x^i} + \sum_jx^j \frac{\partial \varphi}{\partial x^i x^j} = 0.\]
Since $\displaystyle \frac{\partial \varphi}{\partial x^i x^j} = g_{ij}$ we obtain
\[\textrm{grad} \, \varphi = \sum_{i,j} g^{ij} \frac{\partial \varphi}{\partial x^j} \frac{\partial}{\partial x^i} = -\sum_j x^j \frac{\partial}{\partial x^j}.\]
It is equivalent to $(*)$ that
\[\frac{\partial \varphi}{\partial x^j} + \sum_k x^k g_{jk} = 0.\]
Taking the derivative of both sides with respect to $x^i$ and applying Proposition \ref{H} we have
\[g_{ij} + g_{ij} + \sum_k 2x^k \gamma_{ijk} = 0,\]
that is,
\[\sum_k x^k \gamma_{ijk} = -g_{ij,}\]
\begin{flushright}
$\Box$
\end{flushright}
\end{Prf}

We set $\displaystyle H= \sum_j x^j \frac{\partial}{\partial x^j} \, (= -\textrm{grad} \, \varphi)$ and denote by $\mathscr{L}_H$ Lie differentiation with respect to $H$. 

\begin{Prop}
For  $\sigma \in A^p(\Omega)$ we have
\[\mathscr{L}_H \sigma = D_H \sigma + p \sigma.\]
\end{Prop}

\begin{Prf}
For $X \in \mathscr{X}(\Omega) $ we obtain
\[D_H X = X,\]
and so
\[[H,X] = D_H X - D_X H = D_H X - X.\]
Then for $X_1, \dots ,X_p \in \mathscr{X}(\Omega)$ we have
\[
\begin{split}
&(\mathscr{L}_H \sigma)(X_1, \dots ,X_p) \\
=&H \sigma (X_1, \dots ,X_p) - \sum_i \sigma (X_1, \dots , [H,X_i], \dots X_p) \\
=&H \sigma (X_1, \dots ,X_p) - \sum_i \sigma (X_1, \dots ,D_H X_i, \dots ,X_p) + p \sigma (X_1, \dots ,X_p) \\
=&(D_H \sigma) (X_1, \dots ,X_p) + p \sigma (X_1, \dots ,X_p).
\end{split}
\]
\begin{flushright}
$\Box$
\end{flushright}
\end{Prf}

By Cartan's formula we have the following.

\begin{Cor} \label{Lie}
For $\omega \in A^{p,q}(\Omega)$ we have
\[(\partial i(H) + i(H) \partial ) \omega = D_H \omega + p \omega,\]
\[(\bar{\partial} \bar{i}(H) + \bar{i}(H) \bar{\partial} ) \omega = D_H \omega + q \omega.\]
\end{Cor}

\begin{Main} \label{Con}
Let $(\Omega,D,g=Dd\varphi)$ be a regular convex cone in $\mathbb{R}^n$ with the Cheng-Yau metric. Then for $p>q$ and all $v \in L ^{p,q}(M,g)$ such that $ \bar { \partial }v=0$, there exists $u \in L^{p,q-1}(M,g)$ such that 
\[ \bar { \partial }u=v, \quad  \| u \| \le  (p-q)^{- \frac {1}{2}} \| v \| .\]
In particular, we have 
\[L^2H^{p,q}_{\bar{\partial}}(M,g)=0, \quad \textrm{for} \ p>q.\]
\end{Main}

\begin{Prf}
By Theorem \ref{KI} we obtain
\[ \Lambda \partial  + \partial \Lambda = - \bar{\delta} + \bar{i}(X_\alpha), \quad 
\Lambda \bar{\partial} + \bar{\partial} \Lambda = -\delta + i(X_\alpha) .\]
Then we have
\[
\begin{split}
\bar{\partial} \bar{\delta} 
&= \bar{\partial} (- \Lambda \partial -\partial \Lambda + \bar{i}(X_\alpha)) \\
&= (\Lambda \bar{\partial} + \delta - i(X_\alpha)) \partial - \bar{\partial} \partial \Lambda + \bar{\partial} \bar{i}(X_\alpha) \\
&= \delta \partial -i(X_\alpha)\partial + \bar{\partial} \bar{i}(X_\alpha) + \Lambda \bar{\partial} \partial - \bar{\partial} \partial \Lambda, \\
\bar{\delta} \bar{\partial}
&= (-\Lambda \partial - \partial \Lambda + \bar{i}(X_\alpha)) \bar{\partial} \\
&= -\Lambda \partial \bar{\partial} + \partial (\bar{\partial} \Lambda + \delta - i(X_\alpha)) + \bar{i}(X_\alpha) \bar{\partial} \\
&= \partial \delta - \partial i(X_\alpha) + \bar{i}(X_\alpha) \bar{\partial} - \Lambda \partial \bar{\partial} + \partial \bar{\partial} \Lambda,
\end{split}
\]
and so
\[\bar{\square} = \square -(\partial i(X_\alpha) + i(X_\alpha) \partial) +(\bar{\partial} \bar{i}(X_\alpha) + \bar{i}(X_\alpha) \bar{\partial}),\]
where $\square = \partial \delta + \delta \partial$.

Since $\varphi$ is the solution of the equation in Theorem \ref{M},
\[X_\alpha = \textrm{grad} \, \varphi = -H. \]
Hence by Corollary \ref{Lie},
\[\bar{\square} = \square + p - q\]
Therefore, for $\omega \in A^{p,q}_0(\Omega)$ we obtain
\[\| \bar{\partial} \omega \|^2 + \| \bar{\delta} \omega \|^2 \ge (p-q) \| \omega \|^2.\]
Then the assertions are proved similarly with Main theorem \ref{L}.
\begin{flushright}
$\Box$
\end{flushright}
\end{Prf}

We have the following from Main theorem \ref{Con} and Theorem \ref{Reg}.

\begin{Cor} \label{Or}
Let $(\Omega,D,g=Dd\varphi)$ be a regular convex cone in $\mathbb{R}^n$ with the Cheng-Yau metric. Then we have
\[L^2H^{p,q}_{\bar{\partial}}(\Omega)=0, \quad \textrm{for} \ p+q>n \ \textrm{or} \ p>q.\]
\end{Cor}

\subsection{$L^2$-cohomology groups on $\mathbb{R}_+^n$}
Let $(\mathbb{R}_+^n,D,g)$ be the same as in Example \ref{E}. Then we can apply Corollary \ref{Or} to $(\mathbb{R}_+^n,D,g)$. However, we have a stronger vanishing theorem.

\begin{Main} \label{Ep}
For  $p \ge 1$, $q \ge 0$ and $v \in L ^{p,q}(\mathbb{R}_+^n,g)$ such that $\bar { \partial }v=0$, there exists $u \in L^{p,q-1}(\mathbb{R}_+^n,g)$ such that 
\[ \bar { \partial }u=v, \quad  \| u \| \le p^{- \frac {1}{2}} \| v \| .\]
In particular, we have 
\[L^2H^{p,q}_{\bar{\partial}}(\mathbb{R}_+^n,g)=0, \quad \textrm{for} \ p \geq 1 \ \textrm{and} \ q \geq 0.\]
\end{Main}

In this section we show Main theorem \ref{Ep}.
For a canonical coordinate $x=(x^1, \dots ,x^n)$ on $\mathbb{R}_+^n$, we set $t=(t^1, \dots ,t^n)=(\log x^1, \dots, \log x^n)$.

\begin{Lem} \label{t1}
The following equations hold.
\item[$(1)$] $\displaystyle g(\frac{\partial}{\partial t^i}, \frac{\partial}{\partial t^j})=\delta_{ij}$.
\item[$(2)$] $\displaystyle D_\frac{\partial}{\partial t^i} \frac{\partial}{\partial t^j}= \delta_{ij}\frac{\partial}{\partial t^j}, \quad D^*_\frac{\partial}{\partial t^i} \frac{\partial}{\partial t^j}= -\delta_{ij}\frac{\partial}{\partial t^j}$.
\item[$(3)$] $\displaystyle D_\frac{\partial}{\partial t^i} dt^j= -{\delta_i}^j dt^j, \quad D^*_\frac{\partial}{\partial t^i} dt^j= {\delta_i}^j dt^j$.
\item[$(4)$] $\displaystyle \alpha= -\sum_{j} dt^j $, \quad where $\alpha$ is the first Koszul form for $(D,g)$.
\end{Lem}

\begin{Lem} \label{t2}
On $(\mathbb{R}_+^n,D,g)$ we have
\[\bar{\delta} = -\sum_j \mathscr{D}_\frac{\partial}{\partial t^j} \bar{i}(\frac{\partial}{\partial t^j}).\]
\end{Lem}

\begin{Prf}
By Proposition \ref{del}, Lemma \ref{t1} and \ref{scD} we obtain
\[
\begin{split}
\bar{\delta} &= -\sum_j \bar{i}(\frac{\partial}{\partial t^j}) \mathscr{D}_\frac{\partial}{\partial t^j}-\bar{i}(\sum_j \frac{\partial}{\partial t^j}) \\
&= -\sum_j \{\bar{i}(\frac{\partial}{\partial t^j}) \mathscr{D}_\frac{\partial}{\partial t^j}+\bar{i}(D_\frac{\partial}{\partial t^j} \frac{\partial}{\partial t^j})\} \\
&= -\sum_j \mathscr{D}_\frac{\partial}{\partial t^j} \bar{i}(\frac{\partial}{\partial t^j}).
\end{split}
\]
\begin{flushright}
$\Box$
\end{flushright}
\end{Prf}

\begin{Prop} \label{Lap}
We denote $\omega \in A^{p,q}(\mathbb{R}_+^n)$ by $\displaystyle \omega=\sum_{I_p,J_q} \omega_{I_p J_q}dt^{I_p} \otimes \overline{dt^{J_q}}$. Then we have
\[\bar{\square} \, \omega= \sum_{I_p,J_q} (\Delta + p) \, \omega_{I_p J_q}dt^{I_p} \otimes \overline{dt^{J_q}},\]
where $\displaystyle \Delta = -\sum_{j} (\frac{\partial}{\partial t^j})^2 $.
\end{Prop}

\begin{Prf}
It is sufficient to show the equation when $\omega = f \, dt^{I_p} \otimes \overline{dt^{J_q}}$. By Lemma \ref{t1} and \ref{t2} we obtain
\[\bar{\partial} \omega = \sum_{i \in J_{n-q}} \frac{\partial f}{\partial t^i} \, dt^{I_p} \otimes \overline{dt^i} \wedge \overline{dt^{J_q}}-\sum_{i \in I_p \cap J_{n-q}} f \, dt^{I_p} \otimes \overline{dt^i} \wedge \overline{dt^{J_q}},\]
\[
\begin{split}
\bar{\delta} \omega &= -\sum_{j \in J_q} \mathscr{D}_\frac{\partial}{\partial t^j} (f \, dt^{I_p} \otimes \bar{i}(\frac{\partial}{\partial t^j}) \, \overline{dt^{J_q}}) \\
&= -\sum_{j \in J_q} \frac{\partial f}{\partial t^j} \, dt^{I_p} \otimes \bar{i}(\frac{\partial}{\partial t^j}) \, \overline{dt^{J_q}} -\sum_{j \in I_p \cap J_q} f \, dt^{I_p} \otimes \bar{i}(\frac{\partial}{\partial t^j}) \, \overline{dt^{J_q}},
\end{split}
\]
\[
\begin{split}
\bar{\delta} \bar{\partial} \omega 
=&-\sum_{i \in {J_{n-q}}} \sum_{j \in J_q \cup \{i\}} \frac{\partial^2 f}{\partial t^i \partial t^j } \, dt^{I_p} \otimes \bar{i}(\frac{\partial}{\partial t^j}) \, (\overline{dt^i} \wedge \overline{dt^{J_q}}) \\
&-\sum_{i \in J_{n-q}} \sum_{j \in I_p \cap (J_q \cup \{i\})} \frac{\partial f}{\partial t^j } \, dt^{I_p} \otimes \bar{i}(\frac{\partial}{\partial t^j}) \, (\overline{dt^i} \wedge \overline{dt^{J_q}}) \\
&+\sum_{i \in I_p \cap J_{n-q}} \sum_{j \in J_q \cup \{i\}} \frac{\partial f}{\partial t^j } \, dt^{I_p} \otimes \bar{i}(\frac{\partial}{\partial t^j}) \, (\overline{dt^i} \wedge \overline{dt^{J_q}}) \\
&+\sum_{i \in I_p \cap J_{n-q}} \sum_{j \in I_p \cap (J_q \cup \{i\})} f  \, dt^{I_p} \otimes \bar{i}(\frac{\partial}{\partial t^j}) \, (\overline{dt^i} \wedge \overline{dt^{J_q}}),
\end{split}
\]
\[
\begin{split}
\bar{\partial}\bar{\delta}\omega
=&-\sum_{j \in J_q} \sum_{i \in J_{n-q} \cup \{j\}} \frac{\partial^2 f}{\partial t^i \partial t^j } \, dt^{I_p} \otimes \overline{dt^i} \wedge \bar{i}(\frac{\partial}{\partial t^j}) \, \overline{dt^{J_q}} \\
&+\sum_{j \in J_q} \sum_{i \in I_p \cap (J_{n-q} \cup \{j\})} \frac{\partial f}{\partial t^j } \, dt^{I_p} \otimes \overline{dt^i} \wedge \bar{i}(\frac{\partial}{\partial t^j}) \, \overline{dt^{J_q}} \\
&-\sum_{j \in I_p \cap J_q} \sum_{i \in J_{n-q} \cup \{j\}} \frac{\partial f}{\partial t^i} \, dt^{I_p} \otimes \overline{dt^i} \wedge \bar{i}(\frac{\partial}{\partial t^j}) \, \overline{dt^{J_q}} \\
&+\sum_{j \in I_p \cap J_q} \sum_{i \in I_p \cap (J_{n-q} \cup \{j\}}) f \, dt^{I_p} \otimes \overline{dt^i} \wedge \bar{i}(\frac{\partial}{\partial t^j}) \, \overline{dt^{J_q}}.
\end{split}
\]
We denote by $(\bar{\delta} \bar{\partial} \omega)_k$ and $(\bar{\partial}\bar{\delta}\omega)_k$ the $k$-th terms of $\bar{\delta} \bar{\partial} \omega$ and $\bar{\partial}\bar{\delta}\omega$ respectively, where $k=1,2,3,4$. Then we have
\[
\begin{split}
(\bar{\delta} \bar{\partial} \omega)_1 + (\bar{\partial}\bar{\delta}\omega)_1
&=-\sum_{j=1}^n (\frac{\partial}{\partial t^j})^2 f \, dt^{I_p} \otimes \overline{dt^{J_q}}, \\
(\bar{\delta} \bar{\partial} \omega)_2 + (\bar{\partial}\bar{\delta}\omega)_3
&=-\sum_{j \in I_p} \frac{\partial f}{\partial t^j} \, dt^{I_p} \otimes \overline{dt^{J_q}}, \\
(\bar{\delta} \bar{\partial} \omega)_3 + (\bar{\partial}\bar{\delta}\omega)_2
&=\sum_{j \in I_p} \frac{\partial f}{\partial t^j} \, dt^{I_p} \otimes \overline{dt^{J_q}}, \\
(\bar{\delta} \bar{\partial} \omega)_4 + (\bar{\partial}\bar{\delta}\omega)_4
&=\sum_{j \in I_p} f \, dt^{I_p} \otimes \overline{dt^{J_q}} = pf \, dt^{I_p} \otimes \overline{dt^{J_q}}.
\end{split}
\]
This completes the proof.
\begin{flushright}
$\Box$
\end{flushright}
\end{Prf}

\begin{Cor}
For $\omega \in A_0^{p,q}(\mathbb{R}_+^n)$ we have
\[\| \bar{\partial} \omega \|^2 + \| \bar{\delta} \omega \|^2 \ge p\| \omega \|^2.\]
\end{Cor}

\begin{Prf}
We denote $\omega \in A_0^{p,q}(\mathbb{R}_+^n)$ by $\displaystyle \omega=\sum_{I_p,J_q} \omega_{I_p J_q}dt^{I_p} \otimes \overline{dt^{J_q}}$. By Lemma \ref{t1} and Proposition \ref{Lap} we obtain
\[
\begin{split}
\| \bar{\partial} \omega \|^2 + \| \bar{\delta} \omega \|^2
&=(\bar{\square} \, \omega,\omega) \\
&=\sum_{I_p,J_q}(\Delta \, \omega_{I_p J_q} , \omega_{I_p J_q}) + p\| \omega \|^2 \\
&\ge p\| \omega \|^2.
\end{split}
\]
\begin{flushright}
$\Box$
\end{flushright}
\end{Prf}

Using the above, we have Main theorem \ref{Ep} similarly with Main theorem \ref{L}.

\end{document}